\documentclass[11pt]{amsart}
\newtheorem{thm}{Theorem}[section]
\newtheorem{cor}[thm]{Corollary}
\newtheorem{lem}[thm]{Lemma}
\newtheorem{prop}[thm]{Proposition}
\theoremstyle{definition}
\newtheorem{defn}[thm]{Definition}
\theoremstyle{remark}

\numberwithin{equation}{section}
\newcommand{\Real}{\mathbb R}

\newcommand{\esssup}{\operatorname{ess\, sup}}
\newcommand{\loc}{\rm{loc}}
\begin{document}

\title[Compact embedding of weighted Sobolev spaces]{Some necessary and some sufficient
conditions for the compactness
of the embedding of weighted Sobolev spaces
}%
\author{Francesca Antoci}
\address{Dipartimento di Matematica, Politecnico di Torino, C.so Duca degli Abruzzi
24, 10129, Torino}
\email{antoci@calvino.polito.it}

\thanks{I would like to thank Prof. E. Salinelli who suggested me the problem}%

%\date{}%
%\dedicatory{}%
%\commby{}%
% ----------------------------------------------------------------
\begin{abstract}
We give some necessary conditions and sufficient conditions for
the compactness of the embedding of Sobolev spaces
$$W^{1,p}(\Omega,w) \rightarrow L^p(\Omega,w),$$ where $w$ is some
weight on a domain $\Omega \subset \Real^n$.
\end{abstract}
\maketitle
% ----------------------------------------------------------------

\section{Introduction}
In the present paper, we give some necessary and sufficient conditions for
the compactness of the embedding of Sobolev spaces
\begin{equation}\label{emb1} W^{1,p}(\Omega,w)\longrightarrow
L^p(\Omega,w).
\end{equation}
Our investigations were originated by a recent paper appeared in the
Annals of Statistics (\cite{Salinelli}), in which a new definition of
Nonlinear
Principal Components is introduced as follows: if ${\bf X}$ is an
absolutely
continuous random vector on an open connected set $\Omega\subseteq
\Real^n $, with density function $f_{{\bf X}}$, zero expectation
and finite variance, the $j^{th}$ Nonlinear Principal Component of
${\bf X}$ is defined as a solution $\varphi_j$ of the maximization
problem
\begin{equation}\label{max} \max_{\psi({\bf X})\in W^{1,2}(\Omega,
f_{{\bf X}})} \frac{{\bf E}(\psi({\bf X})^2)}{{\bf E}(|\nabla
\psi({\bf X})|^2)} \end{equation} subject to the conditions $$
{\bf E}(\psi({\bf X}),\varphi_s({\bf X}))=0\quad \mbox{for}\;
s=1,2,...,j-1, \,j>1.$$ Moreover, it is required that $\varphi_j$
has zero expectation. \par\ Here ${\bf E}(\psi, \varphi)$ denotes
the usual scalar product in the Hilbert space $L^2(\Omega, f_{{\bf
X}})$ of square integrable functions on $\Omega$ with respect to
the measure $f_{{\bf X}}\,dx$, and $W^{1,2}(\Omega, f_{{\bf X}})$
is the weighted Sobolev space $$ W^{1,2}(\Omega, f_{{\bf
X}})=\left\{u\in L^2(\Omega, f_{{\bf X}})\,|\,\forall i=1,...,n,
\; \partial_i u \in L^2(\Omega, f_{{\bf X}})\right\}.$$ It is well-known
that the compactness of the embedding
\begin{equation}\label{compact}W^{1,2}(\Omega, f_{{\bf X}})\longrightarrow
L^2(\Omega, f_{{\bf X}}) \end{equation} turns out to be essential
in order to prove the existence of an orthonormal set of Nonlinear
Principal Components. \par The problem of the compactness of the
embedding (\ref{emb1}) for weighted Sobolev spaces has been
studied by many authors (see e.g. \cite{Gurka-Opic1},
\cite{Gurka-Opic2}, \cite{Gurka-Opic3}). For a  rich bibliography
on this kind of problems we refer also to \cite{Avantaggiati},
\cite{Kufner}, \cite{Kufner-Opic}, \cite{Turesson}.
\par Nevertheless, the attention has mainly
focused  on those classes of weights which arise in the study of
partial differential equations, such as polynomial weights in
unbounded domains, or, in bounded domains, weights depending on
the distance from the boundary and, as concerns the
Poincar\'e-Wirtinger inequality, weights in $A_p$ classes and
derivatives of quasi-conformal mappings. Under this point of view
an enormous amount of work has been done in the last years and a
complete list of contributions is not possible here. Just to have
an idea, see for example \cite{DKN},\cite{FKS}, \cite{FLW},
\cite{HKM}, \cite{S} and the references therein.\par
\par On the contrary, in the maximization problem (\ref{max}) we
have to deal with general density functions $f_{{\bf X}}$.\par
Moreover, most of the criteria for the compactness of the
embedding (\ref{compact}) which can be found in mathematical
literature are not simple to be handled. In view of possible
applications of Nonlinear Principal Components, it is important,
instead, to have simple, easily applicable necessary and sufficient
conditions on the density function $f_{{\bf X}}$ and on the set
$\Omega$ for the compactness of (\ref{compact}). \par In the
present paper, we give some necessary conditions
and sufficient conditions for the compactness of the embedding
(\ref{emb1}). In section 2, we recall some basic facts
about weighted Sobolev spaces. In section 3, we prove the weighted
versions of some necessary conditions for compactness due to Adams
(\cite{Adams}). In section 4, we prove a sufficient condition for
the compactness of (\ref{emb1}), showing that, under suitable
hypothesis on the weight $w$, if the subgraph $$\Omega_w:= \left\{
(x,y)\in \Real^n \times \Real\,|\,x\in \Omega,\,0<y<w(x) \right\}$$ of the
weight is such that the embedding $W^{1,p}(\Omega_w)\longrightarrow
L^p(\Omega_w)$ is compact, then
the embedding (\ref{emb1}) is compact. In section 5, applying a sufficient
condition for compactness in the non-weighted case due to Adams
(\cite{Adams}), we show some examples to which our sufficient condition is
applicable.
 \section{Preliminary facts}
Let $\Omega$ be an open domain in $\Real^n$. We will denote by
$W(\Omega)$ the set of all real-valued, measurable, a.e. in
$\Omega$ positive and finite functions $w(x)$. Elements in
$W(\Omega)$ are called weight functions. For any
Lebesgue-measurable set $U \subset \Real^n$ and for $w\in
W(\Omega)$, we will denote by $\mu_w(U)$ the Borel measure defined
by $$\mu_w(U)= \int_U w(x)\,dx . $$ As usual, we will denote by
$C^{\infty}_c(\Omega)$ the set of all the smooth, compactly
supported functions in $\Omega$. Moreover, we will denote by
$L^p(\Omega, w)$, for $1\leq p
< \infty$, the set of measurable functions $u=u(x)$ such that
\begin{equation}\label{normap} \|u\|_{L^p(\Omega,w)}=\left( \int_{\Omega}
|u(x)|^p\,w(x)\,dx\right)^{\frac{1}{p}}<+\infty. \end{equation} It is a
well-known fact
that the space $L^p(\Omega, w)$, endowed with the norm
(\ref{normap}) is a Banach space. The following Proposition holds
\begin{prop} Let $1\leq p < \infty$. If the weight $w(x)$ is such that
\begin{equation}\label{L1loc} w(x)^{- \frac{1}{p-1}} \in
L^1_{\loc}(\Omega)\quad \quad \mbox{(in the case $p>1$)}
\end{equation}
\begin{equation}\label{p1} \esssup_{x\in B} \frac{1}{w(x)} < +\infty
\quad \quad \mbox{(in the case p=1)}, \end{equation}
for every ball $B\subset \Omega$, then
\begin{equation}\label{COND} L^p(\Omega,w) \subseteq
L^1_{\loc}(\Omega).\end{equation}
\end{prop}
\par \bigskip
For a proof, see \cite{Kufner-Opic1}, \cite{Turesson}.
\par \bigskip
As a consequence, under condition (\ref{L1loc}) ((\ref{p1})),
convergence in $L^p(\Omega,w)$ implies convergence in
$L^1_{\loc}(\Omega)$. Moreover, every function in $L^p(\Omega, w)$
has distributional derivatives.
\par \bigskip
\begin{defn} Let $\Omega$ be an open set in $\Real^n$, let  $1 \leq p
< \infty$,
and let $w(x)$ be a weight function satisfying condition
(\ref{L1loc}) (resp. (\ref{p1})). We define the Sobolev space
$W^{m,p}(\Omega,w)$ as the set of those functions $u\in
L^p(\Omega,w)$ such that their distributional  derivatives
$D^{\alpha}u$, for $|\alpha|\leq m$, belong to $L^p(\Omega,w)$.
\end{defn}
It is a well-known fact that
\begin{thm} If $w(x)$ fulfills condition (\ref{L1loc})
(resp.(\ref{p1})), $W^{m,p}(\Omega,w)$ is a Banach space.\end{thm}
For a proof, we refer to \cite{Kufner-Opic1}.\par \bigskip
Throughout the paper, we will always assume that $w$ satisfy condition
(\ref{L1loc})
(resp. (\ref{p1})).

 \section{Necessary conditions for compactness}
In this section we derive some necessary conditions for the
compactness of the embedding \begin{equation}\label{emb}
W^{1,p}(\Omega,w) \longrightarrow
L^p(\Omega,w),\end{equation}which are generalizations to the
weighted case of analogous results obtained by Adams
(\cite{Adams}, \cite{Adams2}) for non-weighted Sobolev spaces.\par
Let $w(x)\in W(\Omega)$ be a weight on an open set $\Omega$. Let
$T$ be a tesselation of $\Real^n$ with n-cubes of edge $h$. If
$H\in T$, we will denote by $N(H)$, as in \cite{Adams}, the cube
of side $3h$ concentric with $H$ and having faces parallel to
those of $H$, and by $F(H)$ the fringe of $H$, defined by $$
F(H):=N(H) \setminus H.$$ The natural extension to the weighted case of
the concept of $\lambda$-fatness employed by Adams is given by the
following \begin{defn} Let $\lambda >0$. A cube $H\in T$ is called
$(\lambda,w)$-fat (with respect to $\Omega$) if
\begin{equation}\label{fat} \mu_w(H \cap \Omega) > \lambda \,
\mu_w(F(H)\cap \Omega ).\end{equation} If $H$ is not
$(\lambda,w)$-fat, it is called $(\lambda,w)$-thin.
\end{defn} As in the non-weighted case, the following property holds:
\begin{thm}\label{1} Let $1\leq p < \infty$. If the embedding
(\ref{emb}) is compact, then for every $\lambda
>0$ and for every tesselation $T$ of fixed edge $h$, $T$ has only
finitely many $(\lambda,w)$-fat cubes.\end{thm} \begin{proof} The
thesis follows from an easy extension of the proof of Theorem 6.33
in \cite{Adams}, where the Lebesgue measure $\mu$ is replaced by
$\mu_w$ and $\lambda$-fat cubes are replaced by $(\lambda,w)$-fat
cubes. \par
\end{proof} \noindent {\bf Remark} Theorem \ref{1}
implies that if $w$ is a weight on $\Real^n$ which has the
doubling property, the embedding (\ref{emb}) cannot be compact. In
particular, the embedding (\ref{emb}) for $A_p$ weights on
$\Real^n$ cannot be compact. \par \bigskip \noindent A consequence
of Theorem \ref{1} is the following result, which is the extension
to the weighted case of Theorem 6.37 in \cite{Adams}. The
boundedness of the weight plays a crucial role and cannot be
removed; a counterexample is given by the weight $w(x)=x^{\alpha}$
on $I=(0,1)\subset \Real$ for $\alpha\leq -1$, for which the
embedding is compact (see \cite{Gurka-Opic2}) even if $\int_I
w(t)\,dt =+\infty$.
\begin{thm}\label{2N} If $w(x)$ is bounded and the embedding (\ref{emb}) is
compact, then necessarily
\begin{equation}\label{finitevol} \int_{\Omega} w(x)\,dx
\,<\,+\infty .\end{equation} \end{thm}
\begin{proof} The proof essentially follows the argument in
\cite{Adams}; we will give some details in order to show where the
boundedness of the weight is essential. Let $T$ be a tesselation
of $\Real^n$ by cubes of unitary edge, and let
$\lambda=(2(3^n-1))^{-1}$. If  $P$ is the union of the finitely
many $(\lambda,w)$-fat cubes of $T$, then $\mu_w(P\cap \Omega)\leq
\mu_w(P)<+\infty.$ \par If $H$ is a $(\lambda,w)$-thin cube of
$T$, thanks to the choice of $\lambda$ we can choose $H_1\in F(H)$
such that $\mu_w(H\cap\Omega) \leq \frac{1}{2} \mu_w(H_1 \cap
\Omega).$ Analogously, if also $H_1$ is $(\lambda,w)$-thin, we can
choose $H_2 \in F(H_1)$ such that $\mu_w(H_1 \cap
\Omega)\leq\frac{1}{2} \mu_w(H_2\cap \Omega).$ If going on in this
way we can construct an infinite chain
$\left\{H,H_1,H_2,...\right\}$ of $(\lambda,w)$-thin cubes, then
for every $j\in \mathbb N$ $$\mu_w(H\cap\Omega)\leq
\frac{1}{2^j}\mu_w(H_j \cap \Omega), $$ whence, thanks to the
boundedness of $w$, $$\mu_w(H\cap \Omega) \leq  \frac{C}{2^j}$$
for some positive constant $C$ for every $j\in \mathbb N$; thus,
$\mu_w(H\cap\Omega)=0$. As a consequence, if we denote by
$P_{\infty}$ the union of all the $(\lambda,w)$-thin cubes in $T$
for which it is possible to construct such an infinite chain, $
\mu_w(P_{\infty}\cap \Omega)=0.$ \par Let $P_j$ denote the union
of all $(\lambda,w)$-thin cubes $H\in T$ such that any chain of
this type stops at the j-th step (that is, such that $H_j$ is
$(\lambda,w)$-fat). Following the proof of \cite{Adams}, we get
 $$ \mu_w(P_j \cap \Omega) \leq (2j+1)^n 2^{-j} \mu_w(P\cap\Omega),$$ whence
$$\sum_{j=1}^{+\infty} \mu_w(P_j \cap \Omega) \leq \mu_w(P\cap
\Omega) \sum_{j=1}^{+\infty}(2j+1)^n 2^{-j} < +\infty.$$ Since
$\Real^n =P\cup P_{\infty} \cup P_1 \cup P_2 \cup ...$, the thesis follows.
 \end{proof}\bigskip Some
stronger necessary conditions for compactness are given in the
following Theorem, which is a generalization to the weighted case
of Theorem 6.40 in \cite{Adams}:
\begin{thm}\label{3N} Let $w(x)$ be a continuous, bounded weight.
For every $r> 0$, let $\Omega_r$, $S_r$ be defined as $$
\Omega_r:=\left\{x\in \Omega \,|\,|x|>r \right\}$$
$$S_r:=\left\{x\in \Omega \,|\,|x|=r \right\}.$$ Moreover, let us
denote by $A_r$ the surface area, with respect to the weight $w$,
of $S_r$. Then if the embedding (\ref{emb}) is compact,
\begin{enumerate}
\item for every $\epsilon>0$, $\delta>0$, there exists $R>0$ such
that if $r\geq R$ $$\mu_w(\Omega_r)\leq
\delta \, \mu_w(\left\{x\in \Omega\,|\,r-\epsilon \leq |x| \leq r
\right\});$$
\item if $A_r$ is positive and ultimately decreasing as
$r\rightarrow +\infty$ then for every $\epsilon >0$ $$\lim_{r
\rightarrow + \infty} \frac{A_{r+\epsilon}}{A_r}=0.$$
\end{enumerate}
\end{thm}
\begin{proof} The thesis follows from an easy
 extension of the argument in
\cite{Adams}, where the Lebesgue measure $\mu$ is replaced by
$\mu_w$, $\lambda$-fat cubes are replaced by $(\lambda,w)$-fat
cubes and $A_r$ is computed with respect to the weight $w$. \par
\end{proof}

A consequence of Theorem \ref{3N} is the following Corollary,
whose proof is a simple generalization to the weighted case of the
proof of Corollary 6.41 in \cite{Adams}, and is therefore omitted.
\begin{cor} Let $w(x)$ be a continuous, upper bounded weight. If
the embedding (\ref{emb}) is compact, then for every $k \in
\mathbb Z$ $$\lim_{r \rightarrow + \infty}
e^{kr}\,\mu_w(\Omega_r)=0.$$
\end{cor}\par \bigskip \noindent

\section{A sufficient condition}
Let $w\in W(\Omega)$ be a lower semicontinuous weight defined on
an open set $\Omega\subseteq \Real^n$. Let us suppose that $w$
vanishes only on a closed subset $\Omega_0 \subset \Omega$.
Moreover, let $$ \Omega_{\infty}:=\left\{x\in
\Omega\,|\,w(x)=+\infty \right\}, $$ and suppose that
$\Omega_{\infty}$ is closed. Both $\Omega_0$ and $\Omega_{\infty}$
have (Lebesgue) measure equal to zero. \par Moreover, we suppose
that $w$ is bounded from above and from below by positive
constants on any compact set $K\subset \Omega\setminus (\Omega_0
\cup \Omega_{\infty})$.\par
 Let us denote by $\Omega_w$
the subgraph of the weight $w(x)$, that is, the open set
\begin{equation}\label{sottografico}\Omega_w:=\left\{ (x,y)\in
\Real^n \times \Real \,|\, x\in \Omega, \, 0<y<w(x)\,\right\}
\end{equation}
and  consider the map $$ J:W^{1,p}(\Omega,w) \longrightarrow
W^{1,p}(\Omega_w)$$ defined by $$ (Ju)(x,y)=u(x)\quad
\mbox{a.e.}.$$ $J$ is well-defined. It is not difficult to see
that if $u\in J(W^{1,p}(\Omega,w))$, then the distributional
derivative in the $y$-direction $\nabla_yu$ is equal to zero. $J$
is an isometry of $W^{1,p}(\Omega,w)$ onto $J(W^{1,p}(\Omega,w))$
since for every $u\in W^{1,p}(\Omega,w)$
$$\|Ju\|^p_{W^{1,p}(\Omega_w)}= \int_{\Omega_w}|Ju(x,y)|^p \,dx dy
+ \int_{\Omega_w}|\nabla_x(Ju)(x,y)|^p\,dx\,dy=$$
$$\int_{\Omega}\left(\int_0^{w(x)}|u(x)|^p dy\right)\,dx
+\int_{\Omega}\left(\int_0^{w(x)}|\nabla_xu(x)|^p\,dy\right)\,dx=$$
$$=\int_{\Omega} |u(x)|^p w(x)\,dx + \int_{\Omega}|\nabla_x
u(x)|^p\,w(x)\,dx.$$ We will denote by $W^{1,p}_y(\Omega_w)$ the
set $J(W^{1,p}(\Omega,w))$. Moreover, we will denote by
$L^p_y(\Omega_w)$ the completion of $W^{1,p}_y(\Omega_w)$ with
respect to the norm of $L^p(\Omega_w)$.
\par \bigskip
\begin{lem}If $w(x)$ satisfies the conditions above, then $C^{\infty}_c(\Omega
\setminus (\Omega_0 \cup \Omega_\infty))$ is dense in
$L^p(\Omega,w)$, for $1\leq p<+\infty$.\end{lem}
\begin{proof} It suffices to show that, given $f\in L^p(\Omega,w)$,
for every $\epsilon >0$ there exists $g\in C^{\infty}_c(\Omega
\setminus (\Omega_0 \cup \Omega_\infty))$ such that $\|f-g
\|_{L^p(\Omega,w)}< \epsilon.$ Let  $\left\{\Omega_n\right\}$,
$n\in \mathbb N$, be an exhaustion of $\Omega$, defined by
$$\Omega_n:=\left\{x \in \Omega\,|\, \min \left\{
d(x,\Omega_0),d(x,\Omega_\infty),d(x,\partial
\Omega)\right\}>\frac{1}{n} \right\}.$$ Let $f\in L^p(\Omega,w)$;
for every $\epsilon
>0$, there exists $\overline n$ such that $$
(\int_{\Omega\setminus \Omega_{\overline n}}|f(x)|^p w(x)
\,dx)^{\frac{1}{p}}<\frac{\epsilon}{2}.$$ Since $w$ is bounded
from above and from below by positive constants on $\Omega_{\overline n}$,
$u\in L^p(\Omega_{\overline n},w)$ if and only if $u\in
L^p(\Omega_{\overline n})$ and there exist $C_1,C_2>0$ such that
for every $u\in L^p(\Omega_{\overline n},w)$
$$C_2\|u\|_{L^p(\Omega_{\overline n})}\leq
\|u\|_{L^p(\Omega_{\overline n},w)}\leq C_1
\|u\|_{L^p(\Omega_{\overline n})}.$$ \par Hence
$f_{|\Omega_{\overline n}}\in L^p(\Omega_{\overline n})$. As a
consequence, there exists a function $g\in
C^{\infty}_c(\Omega_{\overline n})\subset C^{\infty}_c(\Omega
\setminus (\Omega_0 \cup \Omega_\infty))$ such that
$\|g-f_{|\Omega_{\overline n}}\|_{ L^p(\Omega_{\overline
n})}<(2C_1)^{-1}\epsilon$. Hence, $\|g-f_{|\Omega_{\overline n}}\|_{
L^p(\Omega_{\overline n},w)}<\frac{\epsilon}{2}$. This implies
$$\|g-f\|_{L^p(\Omega,w)}= \|g-f\|_{L^p(\Omega_{\overline
n},w)}+\|f\|_{L^p(\Omega \setminus\Omega_{\overline n},w)}<\epsilon.$$

\end{proof} Hence, since $C^{\infty}_c(\Omega\setminus (\Omega_0 \cup
\Omega_\infty))\subset W^{1,p}(\Omega,w)$, $W^{1,p}(\Omega,w)$ is
dense in $L^p(\Omega,w)$. \par Since for every $u\in
W^{1,p}(\Omega,w)$ $$\int_{\Omega}|u(x)|^p w(x)\,dx
=\int_{\Omega_w}|Ju(x,y)|^p\,dxdy,$$ we get
\begin{lem} If $w(x)$ satisfies the above conditions, $J$ can be
extended to an isometry
$$\overline{J}:L^p(\Omega,w) \longrightarrow L^p_y(\Omega_w).$$
\end{lem} \par \bigskip
As a consequence:
\begin{thm}\label{subp} Let $w(x)$ be a weight, satisfying the above conditions.
If the subgraph $\Omega_w$ of $w(x)$ is such
that the embedding
$$I_{\Omega_w}:\,W^{1,p}(\Omega_w)\longrightarrow L^p(\Omega_w)$$
is compact, then the embedding (\ref{emb}) is compact.
\end{thm}\begin{proof} The embedding (\ref{emb}) is compact if and only if
the embedding $$I_y: W^{1,p}_y(\Omega_w)\longrightarrow
L^p_y(\Omega_w)$$ is compact. But $I_y$ coincides with $P_y \circ
I_{\Omega_w} \circ I$, where $I$ is the immersion
$$I:W^{1,p}_y(\Omega_w)\longrightarrow W^{1,p}(\Omega_w)$$ and
$P_y$ denotes the ``projection" $$P_y:L^p(\Omega_w)\longrightarrow
L^p_y(\Omega_w),$$ defined for a.e. $x\in \Omega$ by $$(P_y
u)(x):=\frac{1}{w(x)}\int_0^{w(x)}u(x,y)\,dy.$$ Since $P_y$ and
$I$ are continuous, the thesis follows.
\end{proof} We remark that the sufficient condition of Theorem
\ref{subp} is not a necessary condition; as a matter of fact,
whilst the embedding $I_{\Omega_w}$ can be compact only if we
assume a certain regularity for the weight $w$,  the embedding
(\ref{emb}) can be compact even if the weight $w$ is extremely
irregular, once provided that it is controlled from above and
below by a ``regular" weight $\Phi$ . Indeed, the following
Proposition holds
\begin{prop} Let $\Phi$ be a weight for which the
embedding
\begin{equation}\label{comp}W^{1,p}(\Omega,\Phi)\longrightarrow
L^p(\Omega,\Phi)\end{equation} is compact. Let $w(x)$ be a weight
such that there exist $\alpha,\beta>0$ such that a.e. in $\Omega$
$$\alpha \Phi(x)\leq w(x)\leq \beta \Phi(x).$$ Then the embedding
$$W^{1,p}(\Omega,w)\longrightarrow L^p(\Omega,w)$$ is compact.
\end{prop}\noindent \begin{proof} It is immediate that $u\in
L^p(\Omega,\Phi)$ if and only if $u\in L^p(\Omega,w)$ and $$\alpha
\|u\|_{L^p(\Omega,\Phi)}\leq \|u\|_{L^p(\Omega,w)}\leq \beta
\|u\|_{L^p(\Omega,\Phi)}.$$ Analogously, $u\in
W^{1,p}(\Omega,\Phi)$ if and only if $u\in W^{1,p}(\Omega,w)$ and
$$\alpha \|u\|_{W^{1,p}(\Omega,\Phi)}\leq
\|u\|_{W^{1,p}(\Omega,w)}\leq \beta
\|u\|_{W^{1,p}(\Omega,\Phi)}.$$ Hence, if $\left\{u_n\right\}$ is
a bounded sequence in $W^{1,p}(\Omega,w)$, it is bounded also in
$W^{1,p}(\Omega,\Phi)$. Due to the compactness of the embedding
(\ref{comp}), there exists a subsequence $\left\{u_{n_k}\right\}$
such that $u_{n_k}$ converges in $L^p(\Omega,\Phi)$. Hence,
$u_{n_k}$ converges also in $L^p(\Omega,w)$.\end{proof}

\section{Some applications}
 Let us now state some simple applications of Theorem \ref{subp}.
In the non-weighted case the following sufficient condition for
compactness holds :
\begin{thm}\label{Adams}(\cite{Adams})  Let $\Omega$ be an open set in
$\Real^n$.
 If
\begin{enumerate}
\item there exists a sequence
$\left\{\Omega_N^*\right\}_{N=1}^{\infty}$ of open subsets of
$\Omega$ such that $\Omega_N^*\subseteq \Omega_{N+1}^*$ and for
every $N$ the embedding $$W^{1,p}(\Omega_N^*)\longrightarrow
L^p(\Omega_N^*)$$ is compact;
\item there exist a flow $\Phi:U\rightarrow \Omega$ and a constant
$c>0$ such that if
$\Omega_N=\Omega \setminus \Omega_N^*$ then
\begin{enumerate}
\item $\Omega_N \times [0,c]\subset U$ for every $N$;
\item $\Phi_t$ is one-to-one, for every $t$;
\item there exists $M>0$ such that for every $(x,t)\in U$
$$|\partial_t \Phi(x,t)|\leq M;$$
\end{enumerate}
\item the functions $d_N(t)=\sup_{x \in \Omega_N}|\det
J\Phi_t(x)|^{-1}$ satisfy
\begin{enumerate}
\item $\lim_{N\rightarrow \infty}d_N(c)=0$;
\item $\lim_{N\rightarrow \infty}\int_0^cd_N(t)\,dt=0,$
\end{enumerate}
\end{enumerate}
then the embedding $$W^{1,p}(\Omega)\longrightarrow L^p(\Omega)$$
is compact.
\end{thm}
 \par \bigskip We recall that a flow on $\Omega$ is a continuously
differentiable map $\Phi : U \rightarrow \Omega$, where $U$ is an
open set in $\Omega \times \Real$ containing $\Omega \times
\left\{0\right\}$, with $\Phi(x,0)=x$ for every $x \in \Omega$.
Moreover, we denote by $\Phi_t$ the map $$\Phi_t: x \longmapsto
\Phi(x,t),$$ and by $J\Phi_t$ the Jacobian matrix of $\Phi_t$.\par
\bigskip
Theorem \ref{Adams}, together with Theorem \ref{subp}, can be used
to get compactness results for weighted Sobolev spaces. Our first
result is connected with Example 6.49 in \cite{Adams}:
\par
\bigskip
\begin{lem}\label{ex1}
Let $\Omega$ be a bounded domain in $\Real^n$, with $C^{\infty}$
boundary, and let $w(x)\in C^1(\Omega)$ be a weight on $\Omega$,
positive in every compact set $K\subset \Omega$, such that, if we
denote by $r(x)$ the distance $$r(x)=\mbox{dist}(x,
\partial \Omega),$$ near to the boundary $w(x)$ can be expressed as
$$ w(x)= f(r),$$ where $f\in C^1$ is positive, nondecreasing, has
bounded derivative $f'$ and satisfies $\lim_{r \rightarrow
0^+}f(r)=0$; then the embedding (\ref{emb}) is
compact.\end{lem}\begin{proof} Since the boundary is regular,
there exist an open neighbourhood $V$ of $\partial \Omega$ in
$\Omega$, a constant $a>0$ and a diffeomorphism  $$\partial \Omega
\times [0,a) \longrightarrow V,$$ $$ (\xi,r) \longmapsto
x(\xi,r)$$ such that $x(\xi,r) \in
\partial \Omega$ if and only if $r=0$. We can suppose
that $r$ is equal to the distance of $x(\xi,r)$ from the
boundary.\par Since $w$ is strictly positive in $\Omega \setminus V$,
 the embedding (\ref{emb}) is compact if and only if $W^{1,p}(V,w)$
is compactly embedded in $L^p(V,w)$. Let us consider the subgraph
$V_w$ of $w_{|V}$, $$V_w =\left\{(\xi,r,y)\in
\partial \Omega \times [0,a) \times \Real
\,|\,r>0,0<y<f(r)\right\},$$ and, for $N\in \mathbb N$, the sets
$$(V_w)_N:=\left\{(\xi,r,y)\in V_w\,|\,
0<r\leq\frac{1}{N}\right\}.$$ For $N\in \mathbb N$, the open sets
$(V_w)_N^*:=V_w \setminus (V_w)_N$ are such that
$$(V_w)_N^*\subseteq (V_w)_{N+1}^*;$$ moreover, they satisfy the
cone property, hence the embedding
$$W^{1,p}((V_w)_N^*)\longrightarrow L^p((V_w)_N^*)$$ is compact
for every $N \in \mathbb N$. Let us consider the flow
$$\Phi(\xi,r,y,t):= \left(\xi,r+t,\frac{f(r+t)}{f(r)}y\right)$$
defined on the set $$U=\left\{(\xi,r,y,t)\in \partial \Omega\times
[0,a)\times \Real\times \Real\,|\,(\xi,r,y)\in V_w, -r <t<a-r
\right\}.$$ We simply have to check that the sets $(V_w)_N^*$ and
the flow $\Phi$ satisfy the conditions of Theorem \ref{Adams}. It
is easy to see that $(V_w)_N \times [0,\frac{a}{2}] \subset U$ for
every $N\in \mathbb N$. $\Phi_t$ is one-to-one for every $t$,
since $$ \det(J\Phi_t)=\frac{f(r+t)}{f(r)}>0.$$ Moreover,
$$|\partial_t\Phi(\xi,r,y,t)|=|(0,1,\frac{f'(r+t)}{f(r)}y)|\leq
M$$ for some $M>0$ since $f'$ is bounded and $|\frac{y}{f(r)}|<1$
on $V_w$.\par Further, $$d_N(t):=\sup_{(V_w)_N}|\det
(J\Phi_t)|^{-1}=\sup_{(V_w)_N}|\frac{f(r)}{f(r+t)}|$$ satisfies
$$\lim_{N\rightarrow
+\infty}d_N\left(\frac{a}{2}\right)=\lim_{r\rightarrow
0^+}\frac{f(r)}{f(r+\frac{a}{2})}=0.$$ Analogously, for every
$t>0$ $$ \lim_{N\rightarrow +\infty}d_N(t)=0,$$ and by dominated
convergence $$\lim_{N\rightarrow
+\infty}\int_0^{\frac{a}{2}}d_N(t)\,dt =0.$$ Hence, by Theorem
\ref{Adams}, the embedding $I_{V_w}$ is compact, and Theorem
\ref{subp} implies that the embedding (\ref{emb}) is compact.
\end{proof}
{\bf Remark} In particular, for $f(r)=r^{\alpha}$, $\alpha\geq 1$,
we find compact embedding as in \cite{Gurka-Opic2}. This holds
also for $\alpha=p-1$. \par \bigskip \noindent In the case of a
radial weight on $\Real^n$, combining Theorems \ref{Adams},
\ref{subp} and \ref{3N} we can even get a necessary and sufficient
condition. The following result is a sort of ``radial" version of
Example 6.48 in \cite{Adams}:
 \begin{lem}\label{ex2} Let
$\Omega= \Real^n$, and $w(x)$ be a radial function $w(x)=g(r)$
where $r=|x|$ and $g\in C^1([0,+\infty))$ is positive,
nonincreasing, with bounded derivative $g'$; then the embedding
(\ref{emb}) is compact if and only if
\begin{equation}\label{infinity}\lim_{s\rightarrow
+\infty}\frac{g(s+\epsilon)}{g(s)}=0\end{equation} for every
$\epsilon
>0$.
\end{lem}\begin{proof} Suppose, first, that (\ref{infinity})
holds for every $\epsilon >0$. Let us consider, on $\Real^n$,
polar coordinates $(r,\theta)$. Then the subgraph of $w$ can be
described by $$
\Omega_w=\left\{(r,\theta,y)\,|\,0<y<g(r)\right\}.$$ For every
$N\in \mathbb N$, let us consider the set
$$(\Omega_w)_N:=\left\{(r,\theta,y)\in \Omega_w\,|\,r\geq N
\right\}.$$ Then $(\Omega_w)^*_N:=\Omega_w \setminus (\Omega_w)_N$
is bounded and has the cone property; hence, the embedding
$$W^{1,p}((\Omega_w)^*_N)\longrightarrow L^p((\Omega_w)^*_N)$$ is
compact for every $N\in \mathbb N$. Moreover
$(\Omega_w)^*_N\subset (\Omega_w)^*_{N+1}$ for every $N\in \mathbb
N$. \par An easy computation shows that the flow
$$\Phi(r,\theta,y,t):=\left(r-t,\theta,\frac{g(r-t)}{g(r)}y\right)$$
defined on the set $$U:=\left\{(r,\theta,y,t)\,|\,0<t<r\right\},$$
 satisfies the conditions of
Theorem \ref{Adams} (with $c=1$). As a consequence, the embedding
$I_{\Omega_w}$ is compact, and Theorem \ref{subp} yields the
thesis.
\par Conversely, suppose that the embedding (\ref{emb}) is
compact. Then by Theorem \ref{3N}
$$A_r=\int_{|x|=r}g(r)r^{n-1}\,dr d\theta=C(n)r^{n-1}g(r)$$ must
fulfill the condition $$\lim_{r \rightarrow +\infty}
\frac{A_{r+\epsilon}}{A_r}=0.$$ As a consequence, $$\lim_{r
\rightarrow +\infty}\frac{g(r+\epsilon)}{g(r)}=\lim_{r \rightarrow
+\infty}\frac{C(n)(r+\epsilon)^{n-1}g(r+\epsilon)}{C(n)r^{n-1}g(r)}=0.$$
\end{proof}
{\bf Remark} In particular, for $g(r)= e^{\alpha r}$ and for
$g(r)\sim r^\alpha$ as $r \rightarrow +\infty$ ($\alpha<0$) we get
that there is no compactness, as stated in \cite{Gurka-Opic3}.\par
\bigskip \noindent Theorem \ref{Adams}, together with Theorem
\ref{subp}, can be used also to deal with weights which are not
bounded from above. \par
\begin{lem}\label{ex3} Let $\Omega$ be a bounded domain in $\Real^n$, with
$C^{\infty}$ boundary, and let $w(x)\in C^1(\Omega)$ be a weight
on $\Omega$, positive in every compact set $K\subset \Omega$;
moreover, we suppose that, if we denote by $r(x)$ the distance
$$r(x)=\mbox{dist}(x,
\partial \Omega),$$ near to the boundary $w(x)$ can be expressed as
$$ w(x)= f(r),$$ where $f(r)\rightarrow +\infty$ as $r\rightarrow
0^+$, $f$ is strictly decreasing on $0<r<\delta$ for some
$\delta>0$, $|f'(r)|\geq \frac{1}{C}$ and
\begin{equation}\label{f-1} \lim_{y \rightarrow +\infty}
\frac{f^{-1}(y+ \epsilon)}{f^{-1}(y)}=0\end{equation} for every
$\epsilon >0$. Then the embedding (\ref{emb}) is
compact.\end{lem}
\begin{proof} As in Lemma \ref{ex1}, it suffices to show that
$W^{1,p}(V,w)$ is compactly embedded in $L^p(V,w)$, where $V$ is a
tubular neighbourhood of $\partial \Omega$ in $\Omega$. To this
end, consider the subgraph $V_w$ of $w_{|V}$, and, for $N\in
\mathbb N$, the sets $$(V_w)_N:=\left\{ (\xi,r,y)\in V_w\,|\,n\leq
y<f(r)\right\}. $$ The open sets $(V_w)_N^*:=V_w \setminus
(V_w)_N$ and the flow $$\Phi(\xi,r,y,t):=
\left(\xi,\frac{f^{-1}(y-t)}{f^{-1}(y)}r,y-t\right),$$ defined for
$0< t<y<f(r)$, fulfill the conditions of Theorem \ref{Adams}.
Hence, Theorem \ref{subp} implies that the embedding (\ref{emb})
is compact.
\end{proof}
{\bf Remark} The previous Lemma does not apply to the weight
$w(x)= r^\alpha$ when $\alpha <0$. In this case, the compactness
of the embedding (\ref{emb}) has been proved in
\cite{Gurka-Opic2}.\par \bigskip \noindent Via a similar proof it
can be shown that
 \begin{lem}
Let $\Omega$ be an open domain in $\Real^n$ such that $\underline
0 \in \Omega$, and let $w\in C^1(\Omega \setminus
\left\{\underline 0\right\})$ be a weight of the type
$w(x)=f(|x|)$, where $f(s) \rightarrow + \infty$ as $s \rightarrow
0^+$, $f$ is strictly decreasing on $0<s<\delta$ for some
$\delta>0$, $|f'(s)|\geq \frac{1}{C}$ and
\begin{equation}\label{f-12} \lim_{y \rightarrow +\infty}
\frac{f^{-1}(y+ \epsilon)}{f^{-1}(y)}=0\end{equation} for every
$\epsilon >0$.Then the embedding (\ref{emb}) is compact.\end{lem}
Finally, we show an example of weight not belonging to the $A_p$
class on a bounded domain:
\begin{lem} Let us consider the set $\Omega:= (-\frac{1}{2},\frac{1}{2})$
and the weight $$w(x):= \left\{ \begin{array}{ll} 1 \quad &\mbox{if $x
\leq
0$}\\
 (\log \frac{1}{x})^{\frac{1}{2}} \quad &\mbox{if
$x>0$.}\end{array}\right.$$
Then the embedding (\ref{emb}) is compact. \end{lem}
\begin{proof} Consider the subgraph $\Omega_w$. Since the set
$$ \Omega^-_w:= \left\{ (x,y)\in \Real^2 \,|\, x\in
\Omega,x<0,0<y<w(x)\right\}$$  has the cone property, in order to
prove the compactness of the embedding (\ref{emb}), it suffices to
show that the embedding $I_{\Omega^+_w}$ is compact, where $$
\Omega^+_w:= \left\{ (x,y)\in \Real^2\,|\,x\in \Omega, x>0,
0<y<w(x)\right\}.$$ To this purpose, it is not difficult to check
that the subsets $(\Omega^+_w)_N^*:=\Omega^+_w \setminus
(\Omega^+_w)_N$, where $$(\Omega^+_w)_N:= \left\{(x,y)\in
\Omega^+_w \,|\, N<y<w(x)\right\},$$ and the flow $\Phi$, defined
by $$ \Phi(x,y,t):= \left(
\frac{e^{-(y-t)^2}}{e^{-y^2}}x,y-t\right)$$
 for $0<t<y$, satisfy the conditions of Theorem \ref{Adams}.
Via Theorem \ref{subp}, the thesis follows.
\end{proof}
Hence in this case a Poincar\'e-Wirtinger inequality holds even if the
weight does not belong to the $A_p$ class.

% ----------------------------------------------------------------

\end{document}